\documentclass[review]{elsarticle}
\usepackage{graphicx}
\usepackage{hyperref}
\hypersetup{pdfauthor={Name}}
\usepackage{lineno}
\usepackage{float}
\usepackage{multicol,multirow}
\usepackage{multienum}
\usepackage{eucal}
\usepackage[T1]{fontenc}
\usepackage[utf8]{inputenc}
\usepackage[english]{babel}
\usepackage{tikz}
\usetikzlibrary{shapes}
\usetikzlibrary{matrix,arrows,decorations.pathmorphing}
\usepackage{mathtools}
\usepackage{amsmath,amssymb}

\allowdisplaybreaks

\makeatletter
\newcommand{\subjclass}[2][2020]{%
  \let\@oldtitle\@title%
  \gdef\@title{\@oldtitle\footnotetext{#1 \emph{Mathematics subject classification}: #2}}%
}
\newcommand{\keywords}[1]{%
  \let\@@oldtitle\@title%
  \gdef\@title{\@@oldtitle\footnotetext{\emph{Keywords}: #1}}%
}
\makeatother

\newtheorem{theorem}{Theorem}
\newtheorem{lemma}[theorem]{Lemma}
\newdefinition{remark}{Remark}
\newproof{proof}{Proof}
\modulolinenumbers[5]

\journal{Journal of Applied Numerical Mathematics}

\begin{document}

\begin{frontmatter}
\title{Eigenvector centrality and uniform dominant eigenvalue of graph components}
\tnotetext[t1]{This article is a result of the PhD program funded by Makerere-Sida Bilateral Program under Project 316 "Capacity Building in Mathematics and Its Applications".}
%\tnotetext[t2]{The second title footnote which is a longer
%	text matter to fill through the whole text width and
%	overflow into another line in the footnotes area of the
%	first page.}
\author[1]{Collins Anguzu}
\ead{anguzuco@gmail.com}
\author[2]{Christopher Engstr\"{o}m}
\ead{christopher.engstrom@mdh.se}
\author[1]{John Magero Mango \corref{cor1}}
%\author[1]{John Magero Mango \corref{cor1}%
%	\fnref{fn3}}
\ead{mango@cns.mak.ac.ug}
\author[1]{Henry Kasumba}
\ead{kasumba@cns.mak.ac.ug}
\author[2]{Sergei Silvestrov}
\ead{sergei.silvestrov@mdh.se}
\author[3]{Benard Abola}
\ead{benardabola@yahoo.co.uk}

%\cortext[cor1]{John Magero Mango}
%\fntext[fn1]{PhD student in Mathematics, Makerere University.}
%\fntext[fn2]{Senior Lecturer, M{\"a}lardalen University.}
%\fntext[fn3]{Professor of Mathematics, Makerere University - the corresponding author.}
%\fntext[fn4]{Senior Lecturer, Makerere University.}
%\fntext[fn5]{Professor of Mathematics and Applied Mathematics, M{\"a}lardalen University}
%\fntext[fn6]{Lecturer, Gulu University.}
\address[1]{Department of Mathematics, School of Physical Sciences, Makerere University, Box 7062, Kampala, Uganda}
\address[2]{Division of Mathematics and Physics, The School of Education, Culture and Communication, M{\"a}lardalen University, Box 883, 721 23,  V{\"a}ster{\aa}s,Sweden}
\address[3]{Department of Mathematics, Gulu University, Box 166, Gulu, Uganda}
\cortext[cor1]{Corresponding author}
%\title{Eigenvector Centrality and Uniform Dominant Eigenvalue of Graph Components}
%\author{Collins Anguzu\fnref{myfootnote}}
\begin{abstract}
Eigenvector centrality is one of the outstanding measures of central tendency in graph theory. In this paper we consider the problem of calculating eigenvector centrality of graph partitioned into components and how this partitioning can be used. Two cases are considered, firstly where a single component in the graph has the dominant eigenvalue, and secondly when there are at least two components that share the dominant eigenvalue for the graph. In the first case we implement and compare the method to the usual approach (power method) for calculating eigenvector centrality while in the second case with shared dominant eigenvalues we show some theoretical and numerical results.
\end{abstract}
\end{frontmatter}
\linenumbers
\section{Introduction}
\subsection{Ranking vertices in a network}
In the contemporary world of today, to be popular or to have one's products popular or better still to associate with popular figures is every one's desire. The notion of popularity comes along with the idea of hierarchical assembling or ordering of items of certain characteristics. In other words items of a certain particular category are ranked (graded) in order of importance.

One of the first ranking algorithms is the PageRank algorithm. The PageRank algorithm \cite{pasquinelli2009google} was developed by Sergey Brin and Larry Page in 1998 to rank web pages. Since then, the web search engine (now the famous Google) has made rounds in the world of data science.

Alongside PageRank, other measures of central tendency, have emerged over the years. In this paper we will look at another important centrality, namely eigenvector centrality. It has several applications in social networks, capital flow and spread of diseases \cite{bihari2015eigenvector, lu2016vital, spizzirri2011justification, bonacich2001eigenvector}. Various algorithms are used to calculate the eigenvector centrality of a graph, with the easiest and most widely used being the power method \cite{AnderssonSilvestrovAAM2008:mathintsearchengines, jayalatchumy2019novel, ford2014numerical}.

In this paper we derive a method, blended by the usual power iteration for component-wise computation of eigenvector centrality. We then compare this method with the power method of computing eigenvector centrality of the whole graph, in terms of the number of iterations and the computational time. We begin with the case where the strongly connected components have varying dominant eigenvalues. In this case there exists a component that has the overall dominant eigenvalue. Furthermore, we put our attention on components that have the same dominant eigenvalues in a sub-graph that is a single block of strongly connected components. Also of our keen interest is the case where the graph has several independent blocks of strongly connected components, with all components having the same dominant eigenvalue.

It is well known that eigenvector centrality measures the level of importance of a vertex within a graph. This is achieved by awarding scores to every vertex in the graph, which scores are dependent upon the number and strength of connections to the particular vertex.
A vertex joined by links from popular vertices scores highly and hence is popular \cite{AnderssonSilvestrovAAM2008:mathintsearchengines, negre2018eigenvector, newman2008mathematics}. The description given here above qualifies the PageRank as a variant of the eigenvector centrality based on directed graphs since it uses back links.

\subsection{Power Method}
The Power method is an algorithm that defines a sequence of values recursively obtained by the relation
\begin{equation}\label{eqn1}
	\mathbf{x}_{k+1} = \frac{\mathbb{A}^{\top}\mathbf{x}_k}{\|\mathbb{A}^{\top}\mathbf{x}_k\|}.
\end{equation}
The matrix $ \mathbb{A} $ in \eqref{eqn1} is a diagonalisable matrix of dominant eigenvalue $ \lambda $ and $ \mathbf{x} $ is a nonzero vector (called the dominant eigenvector corresponding to $ \lambda $ ), such that
\begin{equation}\label{eqn2}
	\mathbb{A}^{\top}\mathbf{x} = \lambda\mathbf{x}.
\end{equation}
The algorithm in relation \eqref{eqn1} starts with an approximation or a random vector, $ \mathbf{x}_0 $. For the sequence $ \mathbf{x}_k $ to converge, two assumptions must be made:
\begin{itemize}
	\item [(i)] Matrix $\mathbb{A}$ must have a dominant eigenvalue.
	\item [(ii)] The initial vector $ \mathbf{x}_0 $ must have a nonzero component in direction of the dominant eigenvector.
\end{itemize}
To this end, we incline our work to graphs that can be partitioned into components. Details of graph partitioning can however, be obtained from \cite{monien2007approximation, engsilv2019:grapartpagerankchangnetw, schulz2018graph, kim2011genetic, schloegel2000graph, holland1983stochastic, feinberg1980analyzing}. Our research pays attention to eigenvalues of components and computation of eigenvector centrality of  components with different dominant eigenvalues and those with the same dominant eigenvalues. Here we give only a brief summary of partitioning of graphs into diagonal blocks.

\subsection{Graph Partitioning}
In ranking problems, partitioning network structures as strongly connected components plays a crucial role in improving computational drawbacks
\cite{engsilv2019:grapartpagerankchangnetw, engstrom2016pagerank}. In such a case, the corresponding components yield adjacency matrix $\mathbb{A}$ with a block triangular form, that is,
\begin{equation}
	\mathbb{A} =  \begin{pmatrix}
		\rm{W}_{1,1} & \rm{W}_{1,2} & \cdots & \rm{W}_{1,r} \\
		\rm{0}_{2,1} & \rm{W}_{2,2} & \cdots & \rm{W}_{2,r} \\
		\vdots  & \vdots  & \ddots & \vdots  \\
		\rm{0}_{r,1} & \rm{0}_{r,2} & \cdots &\rm{W}_{r,r}
	\end{pmatrix}
\end{equation}
where $\rm{W}_{ii} = n_i$, $\rm i = 1,2, \ldots, \rm r$ and $\rm n_1+ n_2+ \ldots + n_r =n$. Suppose $\rm{W}_{ir}=0,$ $\rm i > \rm r$, then $\mathbb{A}$ is a block triangular matrix. You can get such a block-triangularisation by finding the strongly connected components as described later in Section \ref{sec:imp}.

%\subsection{Notation and abbreviations}
In this section we have highlighted on ranking vertices in a network, the power method and graph partitioning. The rest of the article is organized as follows. In Section \ref{sec:desc} we give a full description of the method we use in this paper, sighting lemmas and theorems involved. We also mention the circumstances under which the method works or fails. Section \ref{sec:imp} is concerned with the implementation of the method with large data, Section \ref{sec:results} is devoted to results of the implementation, and Section \ref{sec:conclusion} contains conclusions of the study.

\section{Description of the method}\label{sec:desc}
\subsection{An overview of the method}
In this section we show the systematic steps in deriving the method to efficiently calculate the eigenvector centrality of a graph following the chain of components. In addition, some numerical examples to illustrate the method are presented.

%We thus begin the section with a few notations that are or shall be used in this paper. The adjacency matrix $\mathbb{A}$ shall have its transposed form denoted by $\mathbb{A}^{\top}$. The symbol $ \lambda_{\max} $ shall denote the maximum of the dominant eigenvalues of matrix $\mathbb{A}$.
Now consider a 2$\times$2 block diagonal matrix $\mathbb{A}$ whose transposed form is \\ $$\mathbb{A}^\top = \left[ \begin{array}{cc} \rm{W}_{1,1}& 0 \\ \rm{W}_{2,1} & \rm{W}_{2,2}  \end{array} \right], $$ where $\rm{W}_{1,1}$ and $\rm{W}_{2,2}$ represent the two block partitions (components) of the network.  When raised to some integer powers $m$ $(\ge 2)$, we get
\begin{align*}
	(\mathbb{A}^\top)^2 &= \left[ \begin{array}{cc} \rm{W}_{1,1}^2& 0 \\ \rm{W}_{2,1}\rm{W}_{1,1}+\rm{W}_{2,1}\rm{W}_{2,2} & \rm{W}_{2,2}^2  \end{array} \right],\\
	(\mathbb{A}^\top)^3 &= \left[ \begin{array}{cc} \rm{W}_{1,1}^3& 0 \\ \rm{W}_{2,1}\rm{W}_{1,1}^2+\rm{W}_{2,2}\rm{W}_{2,1}\rm{W}_{1,1}+\rm{W}_{2,1}\rm{W}_{2,2}^2 & \rm{W}_{2,2}^3  \end{array} \right],\\
	(\mathbb{A}^\top)^4 &= \left[ \begin{array}{cc} \rm{W}_{1,1}^4& 0 \\ \rm{W}_{2,1}\rm{W}_{1,1}^3+\rm{W}_{2,2}\rm{W}_{2,1}\rm{W}_{1,1}^2+\rm{W}_{2,2}^2\rm{W}_{2,1}\rm{W}_{1,1}+\rm{W}_{2,1}\rm{W}_{2,2}^3 & \rm{W}_{2,2}^4  \end{array} \right],
\end{align*}
and so on, while for a 3 $\times$ 3 block matrix,
\begin{align*}
\mathbb{A}^\top & =\left[\begin{array}{ccc} \rm{W}_{1,1}& 0&0 \\ \rm{W}_{21}&\rm{W}_{2,2} & 0\\\rm{W}_{31}&\rm{W}_{32}&\rm{W}_{3,3}  \end{array}\right], \\
(\mathbb{A}^\top)^2 &=\left[ \begin{array}{ccc} \rm{W}_{1,1}^2& 0&0 \\ \rm{W}_{2,1}\rm{W}_{1,1}+\rm{W}_{2,1}\rm{W}_{2,2}&\rm{W}_{2,2}^2 & 0\\\rm{W}_{3,1}\rm{W}_{1,1}+\rm{W}_{3,1}\rm{W}_{3,3}+\rm{W}_{3,1}\rm{W}_{3,2}&\rm{W}_{3,2}\rm{W}_{2,2}+\rm{W}_{3,2}\rm{W}_{3,3}&\rm{W}_{3,3}^2  \end{array} \right]
\end{align*}
and so on.

We now consider a large network that has generated $r$ block components having links within themselves and also between the components except for the lowest component $\rm{W}_{r,r}$, that we assume not to have out links or is a dangling component. 

The transpose of the matrix $\mathbb{A}^{\top}$, for arbitrary $r$ blocks is
\begin{equation}
	\mathbb{A}^{\top} =  \begin{pmatrix}
		\rm{W}_{1,1} & 0 & \cdots & 0 \\
		\rm{W}_{2,1} & \rm{W}_{2,2} & \cdots & 0 \\
		\vdots  & \vdots  & \ddots & \vdots  \\
		\rm{W}_{r,1} & \rm{W}_{r,2} & \cdots &\rm{W}_{r,r}
	\end{pmatrix},
\end{equation} where $\rm{W}_{i,i}$ for $i = 1, \ldots, r$ are square transitions within the components, $\rm{W}_{i,j}$ are transition matrices, (not necessarily square), into the respective diagonal states in the lower levels, $\rm{W}_{i,j} = 0$ for $j > i$. Raising this matrix to an arbitrary integer power $m$ gives

\begin{equation}\label{eqn5+}
	(\mathbb{A}^\top)^{m} =  \begin{pmatrix}
		\rm{W}_{1,1}^m & 0 & \cdots & 0 \\
		\rm{W}_{2,1}^{(m)} & \rm{W}_{2,2}^m & \cdots & 0 \\
		\vdots  & \vdots  & \ddots & \vdots  \\
		\rm{W}_{r,1}^{(m)} & \rm{W}_{r,2}^{(m)} & \cdots &\rm{W}_{r,r}^m
	\end{pmatrix}.
\end{equation}
Considering the case where $m = 2$, we have the simplified form of expression for the new matrix $\rm{W}_{2,1}^{(2)}$ as 
\begin{equation}\label{mequal2}
\rm{W}_{2,1}\rm{W}_{1,1}+\rm{W}_{2,2}\rm{W}_{2,1} = \rm{W}_{2,2}^0\rm{W}_{2,1}\rm{W}_{1,1}^1+\rm{W}_{2,2}^1\rm{W}_{2,1}\rm{W}_{1,1}^0.
\end{equation} 
Similarly, for $m=3$ we have 

$\rm{W}_{3,1}^{(3)} = \rm{W}_{3,1}\rm{W}_{1,1}^2+\rm{W}_{2,2}\rm{W}_{3,1}\rm{W}_{1,1}+\rm{W}_{2,2}^2\rm{W}_{3,1} = \rm{W}_{2,2}^0\rm{W}_{3,1}\rm{W}_{1,1}^2+\rm{W}_{2,2}^1\rm{W}_{3,1}\rm{W}_{1,1}^1+\rm{W}_{2,2}^3\rm{W}_{3,1}\rm{W}_{1,1}^0$.

This summary leads us to the following result.

\begin{lemma}
	Let $\mathbb{A}$ be a 2$\times$2 block diagonal matrix with blocks $W_{1,1}, W_{2,2}$. And let $\rm{W}_{i,j}$,  $i < j$ be transitions from the block $\rm{W}_{1,1}$ to block $\rm{W}_{2,2}$ with $\rm{W}_{i,j}=0$ for $i > j$. Let $\lambda_{\max_1}$ and $\lambda_{\max_2}$ be the dominant eigenvalues of the components $W_{1,1}$ and $ W_{2,2}$ with  $\lambda_{\max_1}$ $>$ $\lambda_{\max_2}$. Then the vector corresponding to the eigenvector centrality of the lower closed component $\rm{W}_{2,2}$ is given by
	\begin{equation}\label{eqn1vecr}
		\vec{r} = \sum_{k=0}^{\infty}\frac{1}{\lambda_{\max_1}}\left(\frac{\rm{W}_{2,2}}{\lambda_{\max_1}}\right)^k\vec{r}_0.
	\end{equation}
\end{lemma}

\begin{proof}
Consider a 2-component graph with the top component $\rm{W}_{1,1}$ having the dominant eigenvalue $\lambda_{\max_1}$. Let the dominant eigenvalue of the bottom component $\rm{W}_{2,2}$ be $\lambda_{\max_2}$ with $\lambda_{\max_1}$ $>$ $\lambda_{\max_2}$. As such, for convergence of the eigenvectors to the maximum of the dominant eigenvalues, we shall divide each matrix of the components by $\lambda_{\max_1}$. Thus the weighted matrix $\mathbb{A}^\top$ raised to power $m$ can be expressed as
$$\mathbb{A}^\top = \left[ \begin{array}{cc} \rm{W}_{1,1}& 0 \\ \rm{W}_{2,1} & \rm{W}_{2,2}  \end{array} \right], \quad
(\mathbb{A}^\top)^m = \left[ \begin{array}{cc} \rm{W}_{1,1}^m& 0 \\ \rm{W}_{2,1}^{(m)} & \rm{W}_{2,2}^m  \end{array} \right]. $$
	Let $\vec{e} = \left[ \begin{array}{c} \vec{e}_1 \\ \vec{e}_2  \end{array} \right] $ denote the eigenvector with all elements equal to one, of the whole network with $\vec{e}_1$ and $\vec{e}_2$ being the normalised eigenvectors of $\rm{W}_{1,1}$ and $\rm{W}_{2,2}$ respectively, corresponding to their dominant eigenvalues. 	
	From the definition of eigenvector centrality, \begin{equation}\label{eqn0} \mathbb{A}^{\top}\vec{e} = \lambda_{\max}\vec{e},\end{equation}
	where $\lambda_{\max}$ is the maximum absolute value of the dominant eigenvalues of the components.
	Consequently,
\begin{eqnarray}\label{eqnp}
	\left(\frac{\mathbb{A}^{\top}}{\lambda_{\max}}\right)^m \vec{e} & = &\nonumber \frac{1}{\lambda_{\max}^m} \left[\begin{array}{cc} \rm{W}_{1,1}^m& 0 \\ \rm{W}_{2,1}^{(m)} & \rm{W}_{2,2}^m  \end{array}\right]\left[\begin{array}{c} \vec{e}_1 \\ \vec{e}_2\end{array}\right],\\
	& = & \frac{1}{\lambda_{\max}^m}\left[\begin{array}{c}{\mathrm{W}_{1,1}^m} \vec{e}_1 \\ {\mathrm{W}_{2,1}^{(m)}} \vec{e}_1 + {\mathrm{W}_{2,2}^m} \vec e_2
	\end{array}\right].
\end{eqnarray}
From the assumption that $\lambda_{\max}=\lambda_{\max_1}$ $>$ $\lambda_{\max_2}$, it follows that the last term in the second row of relation (\ref{eqnp}) reduces to zero. That is, for sufficiently large m, 
$\left(\frac{\rm{W}_{2,2}}{\lambda_{\max_1}}\right)^m\vec{e}_2$ tends to zero.
Hence the relation (\ref{eqnp}) becomes
	\begin{equation*}\lim_{m \rightarrow \infty} \left(\frac{\mathbb{A}^{\top}}{\lambda_{\max_1}}\right)^m\vec{e} = \lim_{m \rightarrow \infty} \frac{1}{\lambda_{\max_1}^m}\left[ \begin{array}{c} {\mathrm{W}_{1,1}^m}\vec{e}_1 \\ {\mathrm{W}_{2,1}^{(m)}}\vec{e}_1  \end{array} \right]. 
	\end{equation*}
Looking at the recurrence relation for ${\mathrm{W}_{2,1}^{(m)}}$ we get
\begin{equation}
\mathrm{W}_{2,1}^{(m)} = \mathrm{W}_{2,1}\mathrm{W}_{1,1}^{m-1} + \mathrm{W}_{2,2}\mathrm{W}_{2,1}^{(m-1)}.	
\end{equation}
Repeating this procedure or following equation (\ref{mequal2}) for arbitrary $m > 2$, gives the following sum:
\begin{equation}
	\mathrm{W}_{2,1}^{(m)} = \sum_{k=0}^{m} \mathrm{W}_{2,2}^k\mathrm{W}_{2,1}\mathrm{W}_{1,1}^{m-k-1}.	
\end{equation}
Since $\vec{e}_1$ is the eigenvector of $\rm{W}_{1,1}$ corresponding to its dominant eigenvalue, we realise that ${\mathrm{W}_{1,1}^{m}}\vec{e}_1$ = $ {\mathrm{W}_{1,1}^{m-1}}{\mathrm{W}_{1,1}}\vec{e}_1$ = $\lambda_{\max_1}\vec{e}_1 $  and so we can write
\begin{equation}\label{eqn2We11}
	\frac{1}{\lambda_{\max_1}^m}{\mathrm{W}_{2,1}^{(m)}}   \vec{e}_1 = \frac{1}{\lambda_{\max_1}^m}\sum_{k=0}^{m}{\mathrm{W}_{2,2}^k}{\mathrm{W}_{2,1}} \lambda_{\max_1} ^{m-k-1} \vec{e}_1 = \sum_{k=0}^{m} \frac{{\mathrm{W}_{2,2}^k}{\mathrm{W}_{2,1}}}{\lambda_{\max_1} ^{k+1}} \vec{e}_1.
\end{equation}
Clearly ${\mathrm{W}_{2,1}}\vec{e}_1$ is vector representing the rank donation to the lower component $\rm{W}_{2,2}$ from $\rm{W}_{1,1}$. We denote the vector ${\mathrm{W}_{2,1}}\vec{e}_1$ in relation \eqref{eqn2We11} by $\vec{r}_0$. Then,
\begin{align}\label{eqn3}
	\frac{1}{\lambda_{\max_1}^m}{\mathrm{W}_{2,1}^{(m)}}\vec{e}_1 = \sum_{k=0}^{m} \frac{{\mathrm{W}_{2,2}^k}}{\lambda_{\max_1} ^{k+1}} \vec{r}_0 = \sum_{k=0}^{m}\frac{1}{\lambda_{\max_1}}\left(\frac{\rm{W}_{2,2}}{\lambda_{\max_1}}\right)^k\vec{r}_0.
\end{align}
If we let $m \rightarrow \infty$ this sum will clearly converge because of our assumption  $\lambda_{\max_1} > \lambda_{\max_2}$ giving the part of the eigenvector belonging to this component for the dominant eigenvalue of the whole graph. Calculating the eigenvector $\vec{r}$ for the lower component as the infinite sum when $m \rightarrow \infty$ completes the proof  
%	For convergence to the vector corresponding to the dominant eigenvalue $\lambda_{\max_1}$ of the top component $\rm{W}_{1,1}$, we divide the result in equation \eqref{eqn3} by $\lambda_{\max_1}$, (See equality \eqref{eqn0}). Thus the eigenvector centrality
%	$\vec{r} $, of the lower component $\rm{W}_{2,2}$ can be written as
%	\begin{equation}\label{eqn5}
%		\vec{r} = \sum_{k=0}^{m}\frac{1}{\lambda_{\max_1}}\left(\frac{\rm{W}_{2,2}}{\lambda_{\max_1}}\right)^k\vec{r}_0,
%	\end{equation}
\begin{equation}\label{eqn1vecre}
	\vec{r} = \sum_{k=0}^{\infty}\frac{1}{\lambda_{\max_1}}\left(\frac{\rm{W}_{2,2}}{\lambda_{\max_1}}\right)^k\vec{r}_0.
\end{equation}
\qed
\end{proof}
The previous result can easily be generalised for $m \geq 2$ and for more than two block components.

At this point we want to show results related to what we have just proved in (\ref{eqn1vecre}).
%The equality \eqref{eqn5+} emphasizes on the Lemma \ref{lemma:basecase}, Theorem \ref{thm:case1} stated here below, without proof. The proof of the lemma and the theorem can be found in our paper \cite{anco}, (submitted).
Let $\mathbb{A}_1$ and $\mathbb{A}_2 $ be matrices of the top and bottom components, respectively. Suppose $\lambda_1$ and $\lambda_2$ are the dominant eigenvalues of the top and bottom components, respectively. Let $\mathbb{B}$ be the link matrix from $\mathbb{A}_1$ to $\mathbb{A}_2$, which is not necessarily a square matrix and the contributing vector $\vec{v} = \mathbb{B}^{\top}\vec{x}_1$ where $\vec{x}_1$ is calculated from the eigenvalue equation $\lambda_1\vec{x}_1 = \mathbb{A}_1^{\top}\vec{x}_1$ and $\vec{v}$ is normalised such that $|\vec{v}| = |\vec{x}_1^{(0)}|$ and $\vec{x}_1^{(0)}$ is the initial iteration value.
The equality  \eqref{eqn5+} emphasizes on the subsequent Lemma \ref{lemma:basecase} and Theorem \ref{thm:case1}.
The proof of these results can be found in our paper \cite{anguzueigenv2022}.

\begin{lemma}\label{lemma:basecase}
	The eigenvector centrality, calculated using non-normalized power iterations, can be obtained by
	\begin{eqnarray}\label{n1m.-}
		\vec{x}_2^{(k)}  = (\mathbb{A}_2^{\top})^k \vec{x}_2^{(0)} + \lambda_1^k \sum_{i=0}^{k-1}\left(\frac{\mathbb{A}_2^{\top}}{\lambda_1}\right)^i\vec{v}.
	\end{eqnarray}
\end{lemma}

\begin{theorem}\label{thm:case1}
	%	The Eigenvector centrality calculated using non-normalized power iterations can be obtained by:
	%	\begin{eqnarray}\label{n1m.-t1}
	%		\vec{x}_2^{(k)}  = (\mathbb{A}_2^{\top})^k \vec{x}_2^{(0)} + \lambda_1^k \sum_{i=0}^{k-1}\left(\frac{\mathbb{A}_2^{\top}}{\lambda_1}\right)^i\vec{v},
	%	\end{eqnarray}
	If $\lambda_1 > \lambda_2$ and the relation in \eqref{n1m.-} holds, then as $k \rightarrow \infty$, the rank vector for the lower component can be expressed as
	\begin{equation}\label{eq:pseries}
		\vec{x}^{(k)}_2 = \sum_{i=0}^{\infty}\left(\frac{\mathbb{A}_2^{\top}}{\lambda_1}\right)^i\vec{v},
	\end{equation}
otherwise the eigenvector of the lower component is expressed as (power method)
	%\end{theorem}
	
	%\begin{theorem}
	%	The Eigenvector centrality calculated using non-normalized power iterations can be obtained by:
	%	\begin{eqnarray}\label{n1mt2}
	%		\vec{x}_2^{(k)}  = (\mathbb{A}_2^{\top})^k \vec{x}_2^{(0)} + \lambda_1^k \sum_{i=0}^{k-1}\left(\frac{\mathbb{A}_2^{\top}}{\lambda_1}\right)^i\vec{v},
	%	\end{eqnarray}
	%	If $\lambda_1 \le \lambda_2$ at the limit $k \rightarrow \infty$, then a normalised rank vector can be written as:
	\begin{equation}
		\vec{x}^{(k)}_2 = \frac{(\mathbb{A}_2^{\top})^k \vec{v}}    {||(\mathbb{A}_2^{\top})^k \vec{v} ||}
	\end{equation} if $\lambda_1 \le \lambda_2$.
\end{theorem}
\subsection{Functionality of the method}
The method described in Section \ref{sec:desc} works for all strongly connected graphs with the condition that either all the vertices have selfloops or none has a selfloop. The method also works for non-intersecting subgraphs in which some components have dominant eigenvalue less than the overall maximum dominant eigenvalue in the graph. The method however, does not work in the case where the graph comprises isolated strongly connected subgraphs or blocks. In such a case the sequential algorithm annuls the ranks of vertices in all other blocks and  maintains only the ranks of the vertices in the lowest component of the last subgraph.
\subsubsection{Disconnected component subgraphs}
{\bf{Theoretical motivation}}
\vspace{0.5cm}

Recall that the matrices considered in this research are weighted adjacency matrices (transition matrices), thus it is natural to introduce the theoretical understanding for Markov processes.
It is well known that communicative states in a Markov process can be interpreted as strongly connected components. Such kinds of components have been studied in good detail by \cite{silvestrov2021perturbed} in the perspective of Markov chains. They computed the stationary distributions of non-intersecting subgraphs using their transition matrices and then merged them based on a probabilistic point of view. Our interest is to determine eigenvector centrality of such a graph with non-intersecting blocks of strongly connected components. In this case we use the adjacency matrices of these isolated blocks and the number of vertices in the graph. The emphasis in this section is on components with the same dominant eigenvalue.
%We can compare the properties of strongly connected components in matrix analysis with communicative states in Markov processes. [ref] clearly elaborate on finding the stationary distribution of such classes as mentioned above with transience and non-transience. From their work, we make use of the conditions $A_{h1}$, $B_1$ and lemmas 4.3.1 and 4.3.2.

%By good comparison we are able to determine the normalised eigenvector centrality of a graph with non-intersecting blocks of components.

\begin{figure}[H]
	\centering
	\begin{multicols}{3}
		
		\begin{tikzpicture}[>=stealth,x=7mm,y=10mm]
			
			\node[shape=circle,draw=black] (1) at (0,0) {1};
			\node[shape=circle,draw=black] (2) at (-1,-2) {2};
			\node[shape=circle,draw=black] (3) at (1,-2) {3};

			\path [->](1) edge node[right] {} (2);
			\path [->](1) edge node[right] {} (3);
		\end{tikzpicture}

		\begin{tikzpicture}[>=stealth,x=7mm,y=10mm]
			\node[shape=circle,draw=black] (4) at (0,2) {4};
			\node[shape=circle,draw=black] (5) at (0,0) {5};

			\path [->](4) edge node[right] {} (5);

		\end{tikzpicture}
		
	\end{multicols}
	\noindent \caption{Two isolated blocks of strongly connected components $\mathcal G^{(1)}, \ldots, \mathcal G^{(5)}$ representing a graph, with all components having the same dominant eigenvalue.} \label{diag2}
\end{figure}
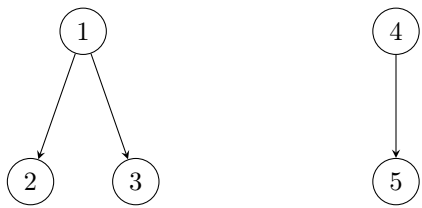

Let $\mathcal G$ be a graph consisting of non-intersecting subgraphs (blocks) $\mathcal G^{(i)}$, $ i = 1, 2, \ldots, r$ with links existing between the components constituting $\mathcal G^{(i)}$. For every $ i = 1, 2, \ldots, r$, the steady state vector $\vec{x}^{(i)}$ is the unique positive solution for the system of linear equation \begin{equation}
	\vec{x}^{(i)} = \mathbb{A}^{\top(i)}\vec{x}^{(i)},
\end{equation} where $\mathbb{A}^{(i)}$ is the adjacency matrix for the block $i$ of components. Let us denote the number of vertices in the subgraph $\mathcal G^{(i)}$ by $k_i$. For $ i = 1, 2, \ldots, r$, we introduce probabilities $(\frac{1}{q})^{(i)}$ such that $\sum_{i=1}^{r}(\frac{1}{q})^{(i)} = \sum_{i=1}^{r}(\frac{k_i}{q}) = 1$, where $q$ is the total number of vertices in the entire graph, $\mathcal G$.  The eigenvector centrality of $\mathcal G$ is then given by \begin{equation}
	\vec{x}^{(i)} = (\frac{1}{q})k_i\vec{x}^{(i)}, ~~  i = 1, 2, \ldots, r.
\end{equation}
For an illustration, we consider a graph with components $\mathcal G^{(1)}, \ldots, \mathcal G^{(5)}$ in isolated subgraphs where all the components have the same dominant eigenvalue. In Figure \ref{diag2} we denote components $\mathcal G^{(1)}$ by 1, $\mathcal G^{(2)}$ by 2, \ldots, $\mathcal G^{(5)}$ by 5.
Let $k_1, k_2, \ldots, k_i$ be the total number of vertices in the components of the $\rm n$ blocks respectively. The overall number of vertices in the entire graph is $q = k_1 + k_2 + \ldots + k_i$. The rank vectors $\vec{x}_1, \vec{x}_2, \ldots, \vec{x}_i$ for the subgraphs are independently obtained. In this case the parallel computation of rank vectors would be the most efficient way. The normalised eigenvector centrality $\vec{x}$ of the entire graph is then approximated by
\begin{equation}\label{eqnEV}
	\vec{x} = \frac{1}{q}\left[\begin{array}{c} k_1\vec{x}_1\\ k_2\vec{x}_2\\ \vdots \\k_i\vec{x}_i \end{array}\right]
\end{equation}
\subsubsection{Example}
Let $\mathbb{A}_1$ and $\mathbb{A}_2$ be the adjacency matrices for the disjointed blocks in Figure \ref{diag2}. Then the matrix $\mathbb{A}$ of the entire graph is of the form
\begin{align*}
\mathbb{A} &=  \left( \begin{array}{cc} \mathbb{A}_1&0\\0&\mathbb{A}_2  \end{array} \right), \\ \mathbb{A}_1 &=  \left( \begin{array}{ccccccc} 1&1&1&1&1&1&0\\1&1&0&1&1&1&1\\0&0&1&1&0&0&0\\0&0&0&1&1&0&0\\0&0&1&0&1&0&0\\0&0&0&0&0&1&1\\0&0&0&0&0&1&1  \end{array} \right), \\
\mathbb{A}_2 &=  \left( \begin{array}{ccccc} 1&1&0&1&1\\0&1&1&1&0\\1&0&1&1&1\\0&0&0&1&1\\0&0&0&1&1  \end{array} \right).
\end{align*}
The corresponding number of vertices in the two blocks are respectively $k_1 = 7$, $k_2 = 5$ giving the total $k = 12$. The eigenvector centrality measures $\vec{x}_1$ and $\vec{x}_2$, of the blocks are independently computed as shown in Table \ref{table:eigv}. In the table also is shown the computed normalised eigenvector centrality $\vec{x}_p$, of the entire graph using relation \eqref{eqnEV}, alongside with the result in the last column $\vec{x}_{np}$ for the unpartitioned graph. It is observed that the two results $\vec{x}_p$ and $\vec{x}_{np}$ compare with an accuracy of about $\pm$ 0.012.
\begin{table}[ht]
	\caption{Eigenvector centrality measure of $\vec{x}_p$ and $\vec{x}_{np}$.} % title of Table
	\centering % used for centering table
	\resizebox{\textwidth}{!}{\begin{tabular}{c c c c c c} % centered columns (4 columns)
		\hline\hline %inserts double horizontal lines
		$\vec{x}_1$ & $\vec{x}_2$ & $k_1\vec{x}_1$ & $k_2\vec{x}_2$ & $\vec{x}_p$&$\vec{x}_{np}$ \\ [0.5ex] % inserts table
		%heading
		\hline % inserts single horizontal line
		$\left( \begin{array}{c} 0\\0\\0.2083\\0.2083\\0.2083\\0.1875\\0.1875  \end{array} \right)$ & $\left( \begin{array}{c} 0\\0\\0\\0.49996\\0.49996  \end{array} \right)$ & $\left( \begin{array}{c} 0\\0\\1.4582\\1.4582\\1.4582\\1.3125\\1.3124  \end{array} \right)$ & $\left( \begin{array}{c} 0.0002\\0.0002\\0.0002\\2.49998\\2.49997  \end{array} \right)$ & $\left(\begin{array}{c} 0\\ 0\\ 0.1282 \\0.1282\\0.1282\\0.1154\\0.1154\\0.1154\\0\\0\\0\\0.1923\\0.1923 \end{array}\right)$&$\left(\begin{array}{c} 0\\ 0\\ 0.1215 \\0.1215\\0.1215\\0.1094\\0.1094\\0.1094\\0\\0\\0\\0.2083\\0.2083 \end{array}\right)$ \\ % inserting body of the table
		[1ex] % [1ex] adds vertical space
		\hline %inserts single line
	\end{tabular}}
	\label{table:eigv} % is used to refer this table in the text
\end{table}

%  \\ $\vec{x}_1 = \left( \begin{array}{c} 0\\0\\0.2083\\0.2083\\0.2083\\0.1875\\0.1875  \end{array} \right)$ and $\vec{x}_2 = \left( \begin{array}{c} 0\\0\\0\\0.49996\\0.49996  \end{array} \right)$. This implies
%$k_1\vec{x}_1 = \left( \begin{array}{c} 0\\0\\1.4582\\1.4582\\1.4582\\1.3125\\1.3124  \end{array} \right)$ and $k_2\vec{x}_2 = \left( \begin{array}{c} 0.0002\\0.0002\\0.0002\\2.49998\\2.49997  \end{array} \right)$. Now the normalised eigenvector centrality is
%$\vec{x} = \frac{1}{q}\left[\begin{array}{c} k_1\vec{x}_1\\ k_2\vec{x}_2\\ \vdots \\k_n\vec{x}_n \end{array}\right] = \left[\begin{array}{c} 0\\ 0\\ 0.1282 \\0.1282\\0.1282\\0.1154\\0.1154\\0.1154\\0\\0\\0\\0.1923\\0.1923 \end{array}\right]$. It compares with $\vec{x} = \left[\begin{array}{c} 0\\ 0\\ 0.1215 \\0.1215\\0.1215\\0.1094\\0.1094\\0.1094\\0\\0\\0\\0.2083\\0.2083 \end{array}\right]$ for computation with the full matrix with an accuracy of $\pm 0.012$.

%Our goal is to find a simple and efficiently calculated expression for the part of the $\vec{x}^{(k+1)} $ vector corresponding to the bottom component depending on the eigenvalue of the top and bottom component when considered alone.

\section{Implementation}\label{sec:imp}

The overall structure of the algorithm can be described in the following steps:
\begin{enumerate}
	\item Component finding: Find the strongly connected components of the graph as well as their levels in the graph. This is done by using a modified version of Tarjan's algorithm as described in \cite{engstrom2016pagerank}.
	\item Intermediate step: Sort components according to level and size for quick access later.
	\item Eigenvector centrality calculation: Eigenvector centrality for each component is calculated individually starting at the top level. Between levels, initial weights for subsequent lower levels are updated. The method is then checked and corrected in case a new largest eigenvalue for a component in the graph is established.
\end{enumerate}

The algorithm is mainly written in Matlab, but with a few selected parts which are written in C or C++ in order to improve the performance of some operations which are comparatively slow in Matlab. The component finding step is done using C++ as a variation of the depth first search implemented in the boost library and parts of the intermediate step is done in C with the rest being implemented in Matlab.

\subsection{Component finding}
For details on the component finding step we refer to \cite{engstrom2016pagerank}. The main difference from Tarjan's algorithm for finding strongly connected components is a little extra bookkeeping (to store levels), as well as updating them at the appropriate steps. We note that our implementation also groups together single vertex components into larger "components" which however, has a very small impact on the results or computation time for the examples considered in this paper. The final complexity of the component finding step is $O(|E|\alpha(|V|))$, in other words amortised linear time. The extra term compared to Tarjan's algorithm for just finding the components comes from the merging of 1-vertex components and properly updating the levels.

\subsection{Intermediate step}
After the strongly connected components of the graph are found the edges and vertices are sorted according to their components level (first) and size(second) or if they lie in between two components and at what level. Note that the sorting into levels here is of linear complexity since we have a known set of (sorted) possible levels so we only have to check how many components we have on each level and assign it to the correct one. This step contains the majority of our overhead in the implementation, but could likely be optimized considerably. Parts of this is written in C due to the comparatively slow performance of some of the operations needed in Matlab.

At this step we also count the number of 1-vertex components on every level so that they can all be handled simultaneously.

\subsection{Eigenvector centrality calculations}
This contains the main part of the algorithm and can be described by the following steps:
\begin{enumerate}
	\item Initialise L as the maximum level among all components and initialize $\lambda_{\max} = -1$ (or some other small number).
	\item For each component of level L: Calculate Eigenvector centrality for that component and dominant eigenvalue $\lambda$ if component did not converge normally.
	\item Check if the largest eigenvalue $\lambda$ among the components of the last level is larger than the previous maximum. If there is a new maximum then zero-out the centrality of all components with a larger lambda as well as initial weights for later components.
	\item Adjust initial weights for all remaining components.
	\item Decrease L by one and go back to step 2 unless we are already at the lowest level in which case we are finished.
\end{enumerate}

We note that the Eigenvector centrality calculation of a component can be adjusted for type and size of the component for more efficient calculation. We have chosen to make just two distinctions here:
\begin{itemize}
	\item Collection of disconnected 1-vertex components: Rather than looping through these one at a time we calculate all of them simultaneously with rank values $C_i = w_i\cdot \lambda_{\max}/(\lambda_{\max}-1)$ if $\lambda_{\max} > 1$ and $C_i=1$ otherwise.
	\item Other components are calculated using the Eigenvector component algorithm.
\end{itemize}

\subsubsection{Eigenvector component algorithm}
The algorithm for calculating Eigenvector centrality for a single component can be described as follows:
\begin{enumerate}
	\item Check if initial weight vector is zero (or sum is below some tolerance).
	\item If step 1 is true, then calculate Eigenvector centrality for the component using power iterations, estimate it's dominant eigenvalue and compare this with $\lambda_{\max}$. If the new eigenvalue is smaller, then set the ranks to zero before exiting.
	\item If the initial vector is not zero then start by iterating using equation \eqref{eqn1vecr}. After a few iterations (we use 20), check if the latest terms in the sum are decreasing (converging) or staying the same or growing (not converging).
	\item If it is determined to converge, let the iterations continue until the maximum change is smaller then the chosen error tolerance.
	\item If it is determined not to converge, break out of the loop and restart the calculations using power iterations instead, making sure the dominant eigenvalue (which should be at least as large as the previous maximum) is stored.
\end{enumerate}

In addition to this we also make two important optimization steps, namely before we start calculating eigenvector centrality in step 2, we compare the largest row-sum with $\lambda_{\max}$, and if this is smaller we can set the ranks to zero and skip the calculations altogether since the resulting dominant eigenvalue has to be smaller than what we already have. Secondly we also check the approximate eigenvalue calculated after some iterations (we use 10 and 20) and again stop if that is lower then half the highest calculated thus far. We use half to be on the safe side, but this is something that would be interesting to investigate further to know how quick you can do the comparison and with what value to do as few iterations as possible.
%\bibliographystyle{plain}
%\bibliography{MyReferences}

\section{Results} \label{sec:results}
All tests are made on a computer with a i5-8250U, 1.6GHz processor. Naturally the properties of the graph will have a large effect on the difference between the baseline algorithm and our componentwise algorithm. The componentwise algorithm gets faster more components, fewer levels (if using a parallel implementation) and if the component with the dominant eigenvalue comes early in the calculation. In particular the best case is if the component with the dominant eigenvalue comes early, is not too large, but has a large enough dominant eigenvalue that other components can immediately be disregarded as having a smaller eigenvalue due to their size.

To test the algorithm we will use a web-graph released by Google as part of a contest 2002 \cite{leskovec2009community}. This graph contains 916428 vertices, 5105039 edges and exhibit many of the properties often found in real world networks such as the presence of a giant component composed of almost half the vertices, as well as exhibiting both a scale free degree distribution and a small-world property.

\begin{figure}[H]
	\begin{center}
		\includegraphics[width=.8\textwidth]{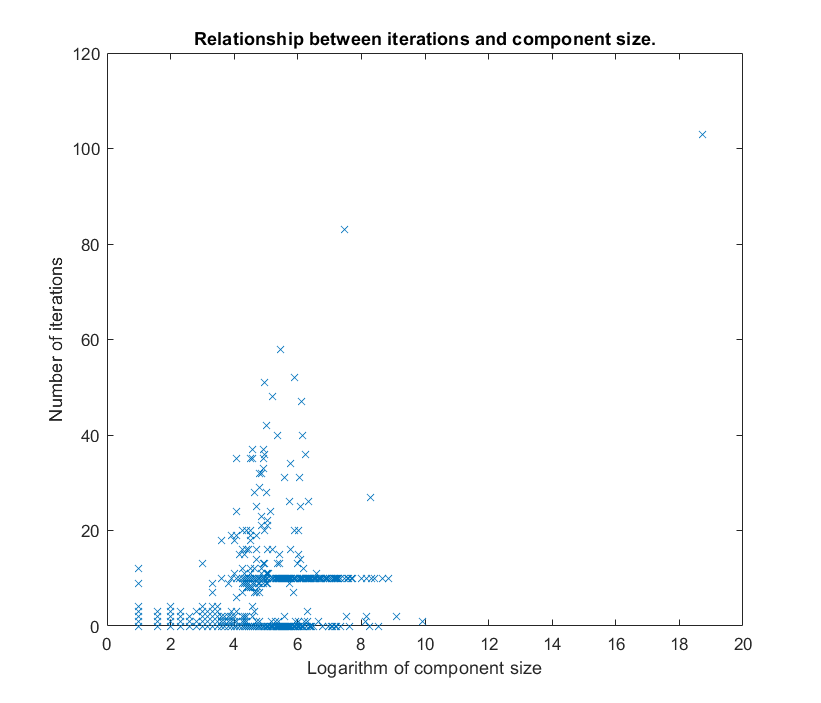}
	\end{center}
	\caption{Number of iterations needed per component plotted against logarithm of the component size.}
	\label{fig:eigen_itr}
\end{figure}

We begin by comparing the number of iterations needed for the baseline algorithm in comparison with the componentwise algorithm. Using an error tolerance of $10^{-9}$, in Figure \ref{fig:eigen_itr} we can see the number of iterations the algorithm needed for each component excluding components of only a single vertex which are handled separately using effectively only 1 iteration. The number of iterations for the components are plotted against the logarithm of the size of the components to also give an overview over how the two relate.

From the figure we can make a few observations:
\begin{itemize}
	\item In comparison with the number of iterations needed for the reference algorithm which was 104, all individual components needed fewer iterations, but the largest component needs just 1 iteration less.
	\item There is a tendency for larger components to need more iterations.
	\item A large number of components stop at 0 and 10 iterations. Those stopping at 0 iterations corresponds to components with zero input vector and a size that is smaller than the dominant eigenvalue thus far. Those stopping at 10 iterations are those where the estimated eigenvalue after 10 iterations are less than half the current max and therefor are discarded. The Figure \ref{fig:eigen_itr} highlights the importance of these checks.
	\item In this graph the giant components are roughly in the middle of the level structure of the graph, if it were calculated earlier it is likely that many other components would be able to be disregarded earlier. In fact, if you know there is a giant component in the graph it may be worthwhile to start the algorithm there since it is almost surely going to be the component with the largest dominant eigenvalue as well.
\end{itemize}

Although the componentwise algorithm requires fewer iterations in components, especially when you also consider 1-vertex components and edges between components (each of which can be considered to do one iteration), the algorithm complexity is much higher since you also need to find components, their levels and reorganize your data. While the number of operations are not too large (each edge is visited twice when finding components and their levels), it is much harder to optimize than the mostly repeated matrix-vector products of the main eigenvector calculation. This result in a significant overhead cost for the componentwise algorithm, the impact of this is not easy to estimate since it so much depends on implementation. It is important to note that our implementation is fairly well optimized, it could likely be optimized further, especially if it was all made in the same programming language with good support for efficient graph search algorithm. On the other hand Matlab is very well optimized for the kind of matrix-vector calculations done in the main eigenvector calculation. On the whole, we expect it to be possible to significantly reduce the relative overhead cost.

\begin{figure}[H]
	\begin{center}
		\includegraphics[width=.8\textwidth]{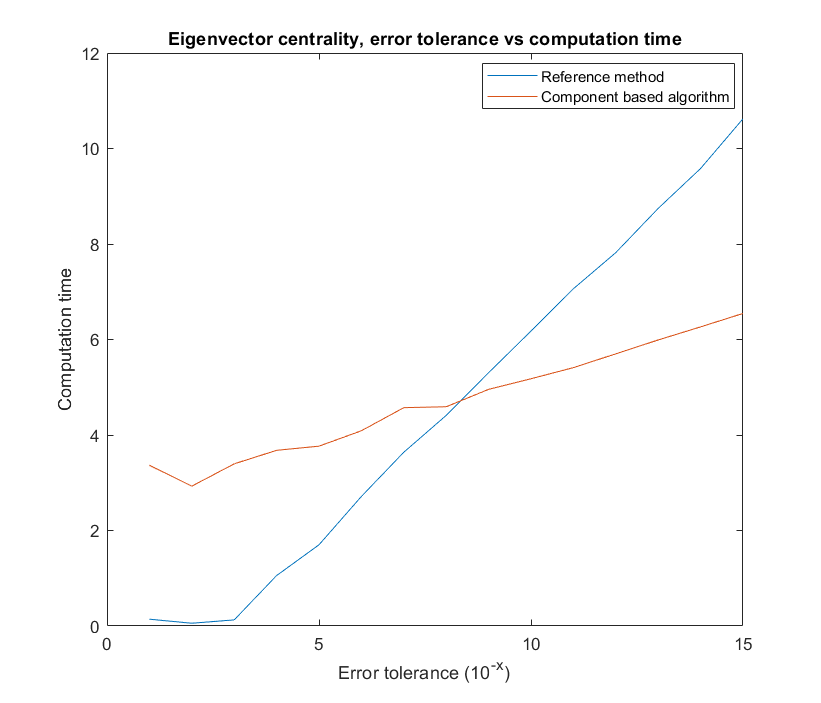}
	\end{center}
	\caption{Computational time versus relative error tolerance for componentwise(red) and reference(blue) method.}
	\label{fig:eigen_tol}
\end{figure}

The overall running time of the two algorithms (componentwise and reference algorithm) for different error tolerances can be seen in Figure \ref{fig:eigen_tol}. As seen here there is a significant overhead for the componentwise algorithm, but the slope is smaller resulting in shorter computation times for small error tolerances. Exactly where the intersection will be obviously depend on the overhead as discussed earlier, but overall the results are positive considering that finding the components also allows for many other possible optimizations we have not done here, but would be possible.

Not shown in the figure here is the fact that, especially for smaller error tolerance, the majority of the time spent calculating the eigenvector centrality of the components is spent in the giant component.

Other optimization that would be possible, but we have not tested here are:
\begin{itemize}
	\item Parallel computation of components on the same level. This has a high potential of increasing the efficiency of the algorithm, especially for graphs without a giant component.
	\item Computing the giant component first. This component is not guaranteed to have the dominant eigenvalue. If it doesn't, then it would result in a significant amount of extra work but when you are  reasonably sure that it will be, this can save a decent amount of computation time. One example where this could be utilized is when you have a persistent network slowly changing over time and you know it had the largest eigenvalue by a significant margin previously.
	\item More efficient handling of specific types of components: We only consider 1-vertex components separately, but it may also be worth to consider other small components, or components with certain properties separately as well to speed up calculation. This is unlikely to make a huge difference to the type of graph we have here since most of the work is already due to the giant component, but may be useful for other types of graphs.
\end{itemize}
\section{Conclusion} \label{sec:conclusion}
In this paper we derive a method, blended by the usual power iteration for componentwise computation of eigenvector centrality. We compare this method with the power method of computing eigenvector centrality of the whole graph (the reference algorithm), in terms of number of iterations and the computational time. We also look at computing eigenvector centrality of strongly connected components of the same dominant eigenvalue and also graphs with independent subgraphs. We found that fewer iterations are needed for the componentwise computation compared to the reference algorithm. Most components do not exceed 10 iterations (Figure \ref{fig:eigen_itr}). The computational time for the two algorithms for varying error tolerances reveals that the componentwise algorithm runs for a shorter time for small error tolerances than the reference algorithm (Figure \ref{fig:eigen_tol}). We have however, found that our method does not work only in some special cases where there are multiple components with the same dominant eigenvalue in disconnected subgraphs. In this case we introduced a probabilistic approach involving the number of vertices in each component.
\section*{Acknowledgment}
This research was supported by the Swedish International Development Corporation Agency (Sida), International Science Program (ISP) in Mathematical Sciences (IPMS), Sida Bilateral Research Program (Makerere University). We also highly recognise the input of the Research environment Mathematics and Applied Mathematics (MAM), Division of Applied Mathematics and Physics, M\"{a}lardalen University, the Department of Mathematics, Makerere University and the Department of Mathematics, Gulu University for providing us with conducive environment for research.

%
%\section{Front matter}
%
%The author names and affiliations could be formatted in two ways:
%\begin{enumerate}[(1)]
%\item Group the authors per affiliation.
%\item Use footnotes to indicate the affiliations.
%\end{enumerate}
%See the front matter of this document for examples. You are recommended to conform your choice to the journal you are submitting to.
%
%\section{Bibliography styles}
%
%There are various bibliography styles available. You can select the style of your choice in the preamble of this document. These styles are Elsevier styles based on standard styles like Harvard and Vancouver. Please use Bib\TeX\ to generate your bibliography and include DOIs whenever available.

%Here are two sample references: \cite{Feynman1963118,Dirac1953888}.

%\bibliographystyle{elsarticle-num}
%\bibliography{johnmangomagero}

\end{document}